\documentclass[a4paper]{article}
\pdfoutput=1
\usepackage{amsmath}
\usepackage{amsthm}
\usepackage{amssymb}
\usepackage{xspace}
\usepackage[utf8]{inputenc}
\usepackage[T1]{fontenc}
\usepackage[l2tabu,orthodox]{nag}
\usepackage{authblk}
\usepackage{verbatim}

\usepackage{hyperref}
\newcommand\rurl[1]{%
  \href{https://#1}{\nolinkurl{#1}}%
}

\pdfpagewidth=\paperwidth
\pdfpageheight=\paperheight
\usepackage[a4paper]{geometry}

\newtheorem{theorem}{Theorem}

\newtheorem{conjecture}[theorem]{Conjecture}

\usepackage[usenames]{color}

\usepackage[mathscr]{eucal}

\title{\vspace{-0.3cm}Counterexamples to a conjecture of Las Vergnas}
\date{}

\author[1]{Robert Brijder}
\author[2]{Hendrik Jan Hoogeboom}
\affil[1]{Hasselt University, Belgium}
\affil[2]{Leiden University, The Netherlands}

\begin{document}

\maketitle

\begin{abstract}
We present counterexamples to a 30-year-old conjecture of Las Vergnas [J.\ Combin.\ Theory Ser.\ B, 1988] regarding the Tutte polynomial of binary matroids.
\end{abstract}

Based on an evaluation established for the Tutte polynomial of plane graphs on $(3,3)$, Michel Las Vergnas made three conjectures in \cite{Vergnas1988367}, in increasing strength, regarding the Tutte polynomial of binary matroids. The first and weakest of these \cite[Conjecture~4.1]{Vergnas1988367} was proved in \cite{Jaeger/Tutte/GFq} and in a more general setting in \cite{TutteMartinOrientingVectors/Bouchet91}.

\begin{theorem}[\cite{Jaeger/Tutte/GFq,TutteMartinOrientingVectors/Bouchet91}]\label{thm:tutte33}
For every binary matroid $M$, the value $T_M(3,3)/T_M(-1,-1)$ is an odd integer.
\end{theorem}

We remark that, for a binary matroid $M$, $T_M(-1,-1) = (-1)^{|E(M)|}(-2)^{b(M)}$, where $b(M)$ is the dimension of the bicycle space of $M$ \cite[Theorem~9.1]{Rosenstiehl1978195}. The other two conjectures remained open for a long time (and were recalled again in 2004 in \cite{EtienneVergnas/MorphismMatroidsIII}). It was shown by Gordon Royle \cite{Royle/VergnasBlog} that $M(K_8)$ is a counterexample to the third and strongest of the three conjectures \cite[Conjecture~4.3]{Vergnas1988367}. In fact, an exhaustive search using the dataset of binary matroids with at most 15 elements of \cite{FripertingerCatalogue} reveals several more counterexamples.


We now state the second conjecture \cite[Conjecture~4.2]{Vergnas1988367}, which is stronger than the first and weaker than the third conjecture. 

\begin{conjecture}[\cite{Vergnas1988367}]\label{conj:LV}
For every binary matroid $M$ and every integer $z$, the value $T_M(-1+4z,-1+4z)/T_M(-1,-1)$ is an odd integer.
\end{conjecture}
Thus Theorem~\ref{thm:tutte33} corresponds to the value $z=1$ in Conjecture~\ref{conj:LV}. It turns out that $M(K_8)$ is not a counterexample to this conjecture. Also, an exhaustive search using the above-mentioned dataset of \cite{FripertingerCatalogue} reveals no counterexample to this conjecture. Consequently, any counterexample has at least 16 elements. In fact, for each binary matroid $M$ with less than 16 elements, $Q_M(z) := T_N(-1+4z,-1+4z)/T_N(-1,-1)$ turns out to have only integer coefficients. This, together with the fact that both $Q_M(0) = 1$ and $Q_M(1)$ are odd (the latter by Theorem~\ref{thm:tutte33}), implies that $Q_M(z)$ is an odd integer for all integers $z$.

Using SageMath~\cite{sagemath} we found that the binary matroid $G$ with 24 elements corresponding to the extended binary Golay code (see, e.g., the appendix of \cite{Oxley/MatroidBook-2nd} for a definition) is a counterexample to Conjecture~\ref{conj:LV}. Moreover, the rank-6 minor $N$ of $G$ with 18 elements having the following reduced representation over $\mathrm{GF}(2)$
\setcounter{MaxMatrixCols}{20}
\[
\begin{pmatrix}
1 & 1 & 0 & 0 & 1 & 0 & 0 & 0 & 1 & 1 & 1 & 1 \\
1 & 0 & 1 & 1 & 0 & 0 & 0 & 1 & 1 & 0 & 1 & 1 \\
0 & 1 & 0 & 0 & 1 & 1 & 1 & 1 & 1 & 1 & 0 & 0 \\
1 & 0 & 1 & 0 & 1 & 1 & 1 & 0 & 1 & 0 & 1 & 0 \\
0 & 1 & 1 & 1 & 0 & 0 & 1 & 0 & 1 & 1 & 1 & 0 \\
0 & 1 & 1 & 0 & 1 & 0 & 1 & 1 & 0 & 0 & 1 & 1
\end{pmatrix}
\]
is another counterexample. Indeed, $N$ has the following Tutte polynomial $T_N(x,y)$
\begin{align*}
y^{12} + 6 y^{11} + 21 y^{10} + 56 y^{9} + 126 y^{8} + 252 y^{7} + x^{6} + 45 x y^{5} + 462 y^{6} + 12 x^{5} + 6 x^{4} y + 225 x y^{4} + \\
747 y^{5} + 72 x^{4} + 111 x^{3} y + 240 x^{2} y^{2} + 675 x y^{3} + 1017 y^{4} + 247 x^{3} + 591 x^{2} y + 1095 x y^{2} + 1057 y^{3} + \\
417 x^{2} + 909 x y + 723 y^{2} + 231 x + 231 y.
\end{align*}
We have $T_N(-1,-1) = 2^6$ and $Q_N(z) = T_N(-1+4z,-1+4z)/T_N(-1,-1)$ is equal to
\begin{align*}
262144 z^{12} - 393216 z^{11} + 344064 z^{10} - 180224 z^{9} + 73728 z^{8} - 18432 z^{7} + 8320 z^{6} - 1248 z^{5} +\\
2616 z^{4} - 1012 z^{3} + \frac{195}{2} z^{2} - \frac{15}{2} z + 1.
\end{align*}
Consequently, $Q_N(z)$ is even for $z \in \{-2, -1, 2\}$, contradicting Conjecture~\ref{conj:LV}.

Finally, the self-dual (but not identically self-dual) rank-9 minor $N'$ of $G$ having the following reduced representation over $\mathrm{GF}(2)$
\[
\begin{pmatrix}
0 & 0 & 0 & 1 & 1 & 1 & 1 & 1 & 1 \\
0 & 1 & 1 & 1 & 0 & 0 & 1 & 1 & 1 \\
0 & 0 & 1 & 0 & 0 & 1 & 0 & 1 & 1 \\
1 & 1 & 0 & 0 & 1 & 0 & 0 & 1 & 1 \\
1 & 1 & 1 & 0 & 0 & 1 & 1 & 1 & 0 \\
1 & 1 & 0 & 1 & 0 & 1 & 0 & 1 & 1 \\
1 & 0 & 1 & 0 & 1 & 0 & 1 & 1 & 1 \\
0 & 1 & 0 & 0 & 1 & 1 & 1 & 1 & 0 \\
1 & 0 & 1 & 1 & 1 & 1 & 0 & 1 & 0
\end{pmatrix}
\]
is yet another counterexample with 18 elements.

\vspace{0.5cm}
\textbf{Acknowledgements}. While SageMath is a collaborative project, we especially would like to thank Stefan van Zwam and Rudi Pendavingh for developing SageMath's matroid theory module, which we used extensively in our calculations. R.B.\ is a postdoctoral fellow of the Research Foundation -- Flanders (FWO).


\bibliographystyle{halpha-abbrv-nodoi}
\bibliography{../mmatroids}

\begin{thebibliography}{ELV04}
\expandafter\ifx\csname url\endcsname\relax
  \def\url#1{\texttt{#1}}\fi
\expandafter\ifx\csname doi\endcsname\relax
  \def\doi#1{\burlalt{doi:#1}{http://dx.doi.org/#1}}\fi
\expandafter\ifx\csname urlprefix\endcsname\relax\def\urlprefix{URL }\fi
\expandafter\ifx\csname href\endcsname\relax
  \def\href#1#2{#2}\fi
\expandafter\ifx\csname burlalt\endcsname\relax
  \def\burlalt#1#2{\href{#2}{#1}}\fi

\bibitem[Bou91]{TutteMartinOrientingVectors/Bouchet91}
A.~Bouchet.
\newblock {Tutte}-{Martin} polynomials and orienting vectors of isotropic
  systems.
\newblock {\em Graphs and Combinatorics}, 7(3):235--252, 1991.

\bibitem[ELV04]{EtienneVergnas/MorphismMatroidsIII}
G.~Etienne and M.~Las~Vergnas.
\newblock The {Tutte} polynomial of a morphism of matroids {III.} {Vectorial}
  matroids.
\newblock {\em Advances in Applied Mathematics}, 32:198--211, 2004.

\bibitem[FW11]{FripertingerCatalogue}
H.~Fripertinger and M.~Wild.
\newblock {\em A catalogue of small regular matroids and their {Tutte}
  polynomials}, 2011.
\newblock \href{https://arxiv.org/abs/1107.1403}{arXiv:1107.1403}. Dataset at
  \rurl{imsc.uni-graz.at/fripertinger/html/matroids/matroide_neu.html}.
  Retrieved: July 24, 2018.

\bibitem[Jae89]{Jaeger/Tutte/GFq}
F.~Jaeger.
\newblock On {Tutte} polynomials of matroids representable over
  $\mathrm{GF}(q)$.
\newblock {\em European Journal of Combinatorics}, 10:247--255, 1989.

\bibitem[LV88]{Vergnas1988367}
M.~Las~Vergnas.
\newblock On the evaluation at $(3,3)$ of the {Tutte} polynomial of a graph.
\newblock {\em Journal of Combinatorial Theory, Series B}, 45(3):367--372,
  1988.

\bibitem[Oxl11]{Oxley/MatroidBook-2nd}
J.~Oxley.
\newblock {\em Matroid theory, Second Edition}.
\newblock Oxford University Press, 2011.

\bibitem[Roy13]{Royle/VergnasBlog}
G.~Royle.
\newblock A {Las} {Vergnas} conjecture, 2013.
\newblock URL
  \rurl{symomega.wordpress.com/2013/06/30/a-las-vergnas-conjecture/}.

\bibitem[RR78]{Rosenstiehl1978195}
P.~Rosenstiehl and R.~Read.
\newblock On the principal edge tripartition of a graph.
\newblock In B.~Bollob\'{a}s, editor, {\em Advances in Graph Theory}, volume~3
  of {\em Annals of Discrete Mathematics}, pages 195--226. Elsevier, 1978.

\bibitem[Sage]{sagemath}
{\em SageMath, the Sage Mathematics Software System}.
\newblock URL \rurl{sagemath.org}.

\end{thebibliography}




\end{document}